\documentclass[12pt,secthm,seceqn]{arXivelsart}

\usepackage{amsmath,amssymb,epsfig}

\newcommand{\Z}{\mathbb{Z}}

\newenvironment{romanlist}
  {\renewcommand{\theenumi}{\roman{enumi}}%
   \setlength{\topsep}{0pt}%
   \vspace{-\parskip}%
   \begin{enumerate}%
     \setlength{\parsep}{0pt}
     \setlength{\parskip}{0pt}%
  }%
  {\end{enumerate}%
   \vspace{-\parskip}}

  {\end{list}}

\DeclareMathOperator{\ord}{ord}

\DeclareMathOperator{\MSD}{MSD}   
\DeclareMathOperator{\lcm}{lcm}

\begin{document}

\begin{frontmatter}
  \title{Noncanonical number systems in the integers}
  \author{Christiaan van de Woestijne}
  \address{Institut f\"ur Mathematik B,
        Technische Universit\"at Graz, 8010 Graz, Austria}

  \begin{abstract}
     The well known binary and decimal representations of the integers, and
     other similar number systems, admit many generalisations. Here, we
     investigate whether still every integer could have a finite expansion on a
     given integer base $b$, when we choose a digit set that does not contain
     $0$. We prove that such digit sets exist and we provide infinitely many
     examples for every base $b$ with $|b|\ge 4$, and for $b=-2$. For the
     special case $b=-2$, we give a full characterisation of all valid digit
     sets.
  \end{abstract}

\end{frontmatter}

\section{Introduction and results}

A \emph{number system} is a coherent notation system for numbers. There are
many possibilities to define such systems, but in this paper we will consider
only generalisations of the \emph{positional} number systems, like the binary
and decimal notations. In such systems, one represents numbers by finite
expansions of the form
\begin{equation} \label{EqFiniteExp}
  a = \sum_{i=0}^\ell d_i b^i,
\end{equation}
where the $d_i$ are taken from a finite set of digits, and $b$ is the
\emph{base} of the system. For example, taking for $b$ an integer greater than
$1$ and using digits $\{0,1,\ldots,b-1\}$, we can represent all nonnegative
integers in the form \eqref{EqFiniteExp}, and these representations are in fact
unique.  However, if we want to represent \emph{all integers} in this form, we
must change either the base or the digit set; for example, we can take an
integer base $b$ with $b\le -2$, and digits $\{0,1,\ldots,|b|-1\}$, as proved
already by Gr\"unwald in 1885 \cite{Grunwald}.

In this paper, we will restrict ourselves to number systems within the set of
integers. The basic definitions are then as follows.

\begin{defn} \label{Def1}
A \emph{pre-number system} in the ring of integers $\Z$ is given by an integer
$b$ and a finite set of integers ${\cal D}$ satisfying the following
properties:
\begin{romanlist}
  \item we have $|b|\ge 2$;
  \item the elements of ${\cal D}$ cover all the cosets of integers modulo 
        $b$.
\end{romanlist}
The integer $b$ is called the \emph{base} of the pre-number system, and
${\cal D}$ is the \emph{digit set}. If $|{\cal D}|=|b|$, we say that
${\cal D}$ is \emph{irredundant}, otherwise it is \emph{redundant}. In an
irredundant digit set, the unique digit that represents the coset of $0$ is
called the \emph{zero digit}.

A pre-number system $(\Z,b,{\cal D})$ is a \emph{number system} if every
$a\in \Z$ has a \emph{finite} expansion of the form
\begin{displaymath} {\textstyle
  a = \sum_{i=0}^{\ell-1} d_i b^i
   }
\end{displaymath}
where all $d_i$ are in ${\cal D}$ and where $\ell$ is a positive integer.

If $(\Z,b,{\cal D})$ is a number system, we call ${\cal D}$ a \emph{valid digit
set} for $b$.
\end{defn}

The notation $(\Z,b,{\cal D})$ for a pre-number system in $\Z$ is motivated
by the fact that pre-number systems may be defined in much more general rings
and other sets (see the forthcoming paper \cite{CvdW2008b}), where instead of
$\Z$ we indicate the set of numbers, or number-like elements, that we want to
have a finite representation. In the present paper, however, all pre-number
systems will be in $\Z$.

Many generalisations of this definition are possible. Already Knuth
\cite[Section 4.1]{Knuth3rd} gave many interesting variants. For all variants
where the basis remains integral in some sense, such as an algebraic integer or
an integer matrix, we would like to refer to Section 3 of the survey paper
\cite{BBLT2006}. It is possible to consider nonintegral bases; this was done in
\cite{Gilbert}, \cite[Section 5.3.3]{WoePhD}, and \cite{AFS2008}. One could
take a positive $b$ and nonnegative digits, and look only at the property of
representing all \emph{nonnegative} integers in the form \eqref{EqFiniteExp};
here, a complete classification of all possible digit sets (which must contain
$0$) was achieved in \cite{Odlyzko1978}, and generalisations to the
higher-dimensional case are given in \cite{LaWa1996-4} and \cite{Akiyama2006}.
There are interesting number systems that use redundant digit sets, such as
those discussed in \cite{Michalek2001,MuSt2005}; in the guise of addition
chains, several such systems are useful for speeding up operations in elliptic
curve cryptosystems (see \cite[Chapter 9]{HandbookElliptic}).

Virtually all papers dealing with number systems as defined above, or with
their generalisations, have used the additional requirement that $0$ be in the
digit set. The main goal of this paper is to explore the consequences when we
drop this restriction, while remaining within the framework of Definition
\ref{Def1}. We will discuss higher-dimensional generalisations in another
paper \cite{CvdW2008b}. Number systems without zero in the case where the base
$b$ is a power of $\pm 2$ were proposed by M\"oller for the purpose of avoiding
Side Channel Attacks in elliptic curve cryptography (see \cite[Section
4.4]{Theriault} and \cite[Section 29.1.1.a]{HandbookElliptic}).

The basic implications of Definition \ref{Def1} will be discussed in Section 2.
For example, if $0$ is not a digit, we cannot pad expansions with zeros if we
want to make them longer; we will be forced to use repetitions of some sequence
of nonzero digits that nonetheless has zero value. We will show that such a
sequence always exists, whenever we have a \emph{number system}. We also show
that the length of such sequences goes to $\infty$ with the size of the
\emph{zero digit}. Next, we construct a few basic examples of digit sets with
without $0$ for any base $b$. Finally, we show that a valid digit set cannot
be translated over an arbitrarily large integer without losing the number
system property, even if it contains $0$ and we leave the $0$ in place. 

In Section 3, we will prove the existence of infinitely many distinct 
sets of nonzero digits in $\Z$ for any integer base $b$ with $|b| \ge 4$, the
main results being Theorems \ref{ThmAvoid} and \ref{ThmAvoidNeg}. This
complements known results for digit sets that do have $0$, which have been
obtained by Matula \cite{Matula} and Kov\'acs and Peth\H o \cite{KoPe1983}. 

As for bases with $|b|\le 3$, we have a pre-number system if $b=\pm 2$ or
$b=\pm 3$. Now for $b=2$, no digit set at all will yield a number system,
whether including $0$ or not; see Corollary \ref{CorOneCycle} for a proof.
For $b=-2$, in Section 4 we will characterise \emph{all} possible digit sets
that yield a number system in $\Z$; although infinite in number, it will turn
out that their structure is different from the infinite families obtained for
larger bases in Section 3. The main result is Theorem \ref{ThmMinus2}. For
$|b|=3$, we have been unable to obtain the existence of infinitely many digit
sets without zero, which therefore remains an open problem.

\section{Digit sets with and without zero}

We will now explore the consequences of not having $0$ as a digit in a number
system. First, we extend some well known results and definitions to the more
general context defined above; see \cite[Sections 2.1, 2.2, 3.1, and
3.2]{BBLT2006} and references therein for more background on these notions.

\subsection{Notations and extensions} \label{SecNotation}

Let $(\Z,b,{\cal D})$ be a pre-number system. For the rest of the paper, we
will assume that all digit sets are \emph{irredundant}. It follows that, given
$a\in\Z$, there exists a unique digit $d_a\in{\cal D}$ such that $a-d_a$ is
divisible by $b$. 

In particular, there will be a unique digit that is itself divisible by $b$;
this is the digit corresponding to the integer $0$, and, as in Definition
\ref{Def1}, we will call it the \emph{zero digit}, whether it be equal to $0$
or not.

\begin{defn} \label{DefT}
Given a pre-number system $(\Z,d,{\cal D})$, define maps
\begin{equation}
  \begin{aligned}
  d &: \Z\rightarrow {\cal D} : a \mapsto d \in{\cal D} \text{ such that }
    b \text{ divides } a - d; \\
  T &: \Z\rightarrow \Z : a \mapsto (a-d(a))/b.
  \end{aligned}
\end{equation}
\end{defn}

The map $T$ is called the \emph{dynamic mapping} of $(\Z,b,{\cal D})$. The name
obviously comes from dynamical system theory; this connection is given in more
detail in \cite{AkiThu2004}. The \emph{digit function} $d$ can also be viewed
as a redefinition of the usual modulo operator: we could say that $d(a)$ is $a$
modulo $b$, with respect to the digits ${\cal D}$.

We will sometimes use the notation $\, a \rightarrow a' \, $ whenever we have
$T(a)=a'$.

\begin{thm} \label{ThmChar}
  A pre-number system $(\Z,b,{\cal D})$, with dynamic mapping $T$, is a number
  system if and only if, for all $a\in \Z$, we have $T^i(a)=0$ for some $i\ge
  1$.
\end{thm} 

\begin{pf*}{Proof.}
For any $a\in \Z$, we want to find the expansion 
\begin{equation} \label{EqExpansion2}
  {\textstyle a = \sum_{i=0}^{\ell-1} d_i b^i }
\end{equation}
with digits in ${\cal D}$ and $\ell\ge 1$. Now the proof is easily done by
induction on $\ell$.
\qed \end{pf*}

The considerations just given show that whether a given pre-number system has
the number system property depends on the structure of the \emph{discrete
dynamical system} on $\Z$ given by the map $T$.

The characterisation given in Theorem \ref{ThmChar} can be made into a finite
algorithm for deciding the number system property, because the dynamical system
just defined is contractive and therefore has a finite \emph{attractor set}
${\cal A}$ \cite{AkiThu2004}. The set ${\cal A}$ by definition has the property
that for all $a\in\Z$, we have $T^n(a)\in{\cal A}$ for $n$ sufficiently large,
and also that $a\in{\cal A}$ implies $T(a)\in{\cal A}$.

Now because the attractor ${\cal A}$ is a finite set, the sequence
$(T^i(a))_{i\ge 0}$ must be purely periodic for any $a\in{\cal A}$; the
elements of ${\cal A}$ that constitute one full period are called a
\emph{cycle} in ${\cal A}$. In the notation given at the beginning of the
section, we can write a cycle in ${\cal A}$ as
\begin{displaymath}
  a_0 \rightarrow a_1 \rightarrow \ldots \rightarrow a_n = a_0,
\end{displaymath}
where $a_{i+1}=T(a_i)$ for all $i$.

The following Theorem is the extension, to general digit sets, of the usual
formulation that in a number system the attractor should contain just
the element $0$ (given as Theorem 3 in \cite{KoPe1983}).

\begin{thm} \label{ThmAlg}
  The pre-number system $(\Z,b,{\cal D})$ is a number system if and only if
  the attractor ${\cal A}$ consists of exactly one cycle under the map $T$,
  and this cycle contains $0$.
\end{thm}

\begin{pf*}{Proof.}
  We have seen that $a\in \Z$ has a finite expansion if and only if 
  $T^i(a)=0$ for some $i\ge 1$. Now if $0\not\in{\cal A}$ and $a\in{\cal A}$,
  then $T^i(a)\ne 0$ for all $i\ge 0$, so that $a$ cannot have a finite
  expansion, and if $a$ is contained in some cycle in ${\cal A}$ that does not
  pass through $0$, we also have $T^i(a)\ne 0$ for all $i$.

  Conversely, if $a\in \Z$, then $T^n(a)\in{\cal A}$ whenever $n$ is large
  enough. Thus if the attractor has just one cycle that also contains $0$,
  there must exist some $i\ge 1$ with $T^i(a)=0$, as desired.
\qed \end{pf*}

The Theorem in particular disallows $1$-cycles in the attractor other than
$0\rightarrow 0$. The next Lemma gives a well-known characterisation of such 
cycles, to be used later.

\begin{lem} \label{Lem1Cycle}
  Let $(\Z,b,{\cal D})$ be a pre-number system, with attractor ${\cal A}$.
  Then ${\cal A}$ contains a $1$-cycle $\, a\rightarrow a\, $ for some $a\in
  \Z$ if and only if $(1-b)a$ is an element of the digit set ${\cal D}$.
\end{lem}

\begin{pf*}{Proof.}
  Let $d\in{\cal D}$, and suppose $d=(1-b)a$ for some $a\in \Z$. It follows
  that
  \begin{displaymath}
    T(a) = (a-d)/b = a,
  \end{displaymath}
  so that ${\cal A}$ has the $1$-cycle $a\rightarrow a$. Conversely, if
  $a\rightarrow a$, then by definition
  \begin{displaymath}
    a = T(a) = (a-d)/b,
  \end{displaymath}
  so we find $d=(1-b)a$.
\qed \end{pf*}

\begin{cor} \label{CorOneCycle}
  Let $(\Z,b,{\cal D})$ be a number system. Then ${\cal D}$ contains no
  nonzero multiples of $1-b$. A fortiori, $|1-b|\ne 1$.
\end{cor}

\begin{pf*}{Proof.}
  Suppose $d=(1-b)a$ for some $a\in \Z$, where $d\in{\cal D}$ is nonzero.
  Then by Lemma \ref{Lem1Cycle},
  \begin{displaymath}
    a\rightarrow a
  \end{displaymath}
  is a nontrivial $1$-cycle in the attractor ${\cal A}$, which contradicts 
  Theorem \ref{ThmAlg}. Furthermore, if $1-b$ is a unit in $\Z$, then obviously
  all digits are multiples of $1-b$, which contradicts the first claim.
\qed \end{pf*}

We are naturally interested in \emph{bounding the size} of the attractor. The
first bound that we will use is well known, and we leave the proof to the
reader.

\begin{lem} \label{LemBd}
  Let $(\Z,b,{\cal D})$ be a pre-number system with dynamic mapping $T$, 
  let $K=\max_{d\in{\cal D}} |d|$, and let $L=K/(|b|-1)$. Let $a\in\Z$.
  \begin{romanlist}
    \item If $|a|>L$, then $|T(a)|<|a|$.
    \item If $|a|\le L$, then also $|T(a)|\le L$.
  \end{romanlist}
\end{lem}

~

Lemma \ref{LemBd} of course implies that $|a|\le L$ for all $a\in{\cal A}$;
however, for pre-number systems in the integers, we can do better than this.
The bounds in Theorem \ref{ThmZBd} below are due to D. Matula \cite[Lemma
6]{Matula} for the case where $0\in{\cal D}$. For the general case, Matula's
argument breaks down, so we will reprove the result. We will use the following
definition, which is interesting in its own right.

\begin{defn} \label{DefNFold}
  Let $(\Z,b,{\cal D})$ be a pre-number system and $n$ a positive integer. We
  define the \emph{$n$-fold digit set} as
  $$
    {\cal D}^n = \left\{ \sum_{i=0}^{n-1} d_ib^i \mid d_i\in{\cal D} \right\}
  $$
  and the \emph{$n$-fold pre-number system} as $(\Z,b^n,{\cal D}^n)$.
\end{defn}

Note that ${\cal D}^n$ is a complete system of representatives of $\Z$ modulo
$b^n$ if and only if ${\cal D}$ is such a system modulo $b$. It follows that
the $n$-fold pre-number system is well defined. The next result gives some
properties of such systems. 

\begin{prop} \label{PropNFold}
  Let $(\Z,b,{\cal D})$ be a pre-number system with dynamic mapping $T$ and
  attractor ${\cal A}$, and let $n$ be a positive integer. Then:
  \begin{romanlist}
    \item
      The dynamic mapping of $(\Z,b^n,{\cal D}^n)$ is equal to $T^n$.
    \item
      The attractor of $(\Z,b^n,{\cal D}^n)$ is equal to ${\cal A}$.
    \item
      $(\Z,b^n,{\cal D}^n)$ is a number system if and only if $(\Z,b,{\cal D})$
      is a number system and $\gcd(n,|{\cal A}|)=1$.
  \end{romanlist}
\end{prop}

\begin{pf*}{Proof.}
  Let $\tilde{T}$ be the dynamic mapping of $(\Z,b^n,{\cal D}^n)$. For all
  $a\in\Z$, we have
  $$
    \tilde{T}(a) = \frac{a - \sum_{i=0}^{n-1} d_ib^i}{b^n},
  $$
  where the digits $d_0,\ldots,d_{n-1}\in{\cal D}$ are chosen so as to make the
  numerator divisible by $b^n$. Thus clearly $\tilde{T}$ is equal to the
  $n$-fold composition of $T$ with itself, as claimed.

  Now let $a\in\Z$ be periodic under $T$ with period length $\ell$; then $a$
  is periodic under $T^n$ with period length
  $\lcm(n,\ell)/n=\ell/\gcd(n,\ell)$. Conversely, if $a$ is periodic under
  $T^n$ with period length $\ell$, then $a$ is also periodic under $T$, with
  some period length that divides $n\cdot \ell$. This proves (ii).
  
  Part (iii) is an easy consequence of (i), together with Theorem \ref{ThmAlg};
  one notes that a cycle of length $\ell$ in ${\cal A}$ is broken up into
  pieces of length $\ell/\gcd(n,\ell)$ if we replace $T$ by $T^n$.
\qed \end{pf*}

\paragraph*{Examples.}
  Theorem \ref{ThmMinus2} implies that $\{1,2\}$ is a valid digit set for the
  base $-2$. The $n$-fold digit set ${\cal D}^n$ is equal to
  $\{-2^n+1,\ldots,-1,0\}$ if $n$ is even and to $\{1,2,\ldots,2^n\}$ if $n$
  is odd. The Proposition now tells us that ${\cal D}^n$ is valid for base
  $(-2)^n$ precisely for odd $n$. In fact, the attractor for all $n$ is equal
  to $\{0,1\}$, but for even $n$ the $2$-cycle $0\rightarrow 1\rightarrow 0$ is
  broken up into two $1$-cycles, and the criterion of Theorem \ref{ThmAlg} is
  violated. One could also have used the obvious criterion that any valid digit
  set for a positive base must contain both negative and positive digits.

  When the starting digit set ${\cal D}$ contains $0$, the attractor ${\cal A}$
  is just $\{0\}$, and the condition on the gcd in (iii) is trivially
  satisfied. Thus, when $0\in {\cal D}$, $(\Z,b,{\cal D})$ is a number system
  if and only if all its $n$-fold pre-number systems are number systems; this
  is Lemma 4 in \cite{Matula}.

\begin{thm} \label{ThmZBd}
  Let $(\Z,b,{\cal D})$ be a pre-number system with attractor ${\cal A}$, and
  let $d=\min_{d\in{\cal D}} d$ and $D=\max_{d\in{\cal D}} d$. Then for all
  $a\in{\cal A}$, we have
  \begin{align*}
       \frac{-D \;\;}{b-1} \;\;  &\le a \le \;\;\; \frac{-d \;\;}{b-1} & 
        \text{ if } b>0; \\[\smallskipamount]
       \frac{-db-D}{b^2-1} &\le a \le \frac{-Db-d}{b^2-1} &
        \text{ if } b<0.
  \end{align*}
\end{thm}

\begin{pf*}{Proof.}
  The proof when $b>0$ is easy and is left to the reader. For the case $b<0$,
  we use the $2$-fold pre-number system $(\Z,b^2,{\cal D}^2)$, which by
  Proposition \ref{PropNFold} has the same attractor as $(\Z,b,{\cal D})$, but
  with the positive base $b^2$. Furthermore, the largest digit of ${\cal D}^2$
  is given by $kb+K$ and the smallest by $Kb+k$, because $b$ is negative. Thus,
  we are reduced to the case of a positive base.
\qed \end{pf*}

\paragraph*{Remark.}
  The interval $[-L,L]$ of Lemma \ref{LemBd} has the property that
  $|a|\le L$ implies $|T(a)|\le L$; we will use this property in Lemma
  \ref{LemKovacs} below. The intervals in Theorem \ref{ThmZBd} only have this
  property for $b>0$. If $I$ denotes the interval given in Theorem \ref{ThmZBd}
  for a negative $b$, then $a\in I$ does imply $T^2(a)\in I$, but we may have
  $T(a)\not\in I$.

\subsection{Zero expansions}

If, in any number system, we have a digit $0$ at our disposal, it is clear that
we can extend any finite expansion for $a$ to any length that we like, by
putting zeros in front.  We now prove an analogous property for a number system
with any given digit set, although we will need repeated instances of a
sequence of more than one digit long to obtain the same effect as zero padding.

\begin{defn}
  A \emph{zero expansion} of a pre-number system $(\Z,b,{\cal D})$ is a
  sequence of digits $(d_0,\ldots,d_{\ell-1})$ in ${\cal D}$, with $\ell\ge 1$,
  such that
  \begin{equation} \label{EqZeroCycle}
    {\textstyle  \sum_{i=0}^{\ell-1} d_i b^i = 0.  }
  \end{equation}
\end{defn}

Note that a zero expansion is already determined by its length; in particular,
if a pre-number system has a zero expansion at all, then it also has a 
\emph{shortest} zero expansion, which is uniquely determined.

\begin{thm} \label{ThmZeroExpansion}
  Every number system $(\Z,b,{\cal D})$ has a unique zero expansion of minimal
  length.
\end{thm}

\begin{pf*}{Proof.}
  We use Theorem \ref{ThmAlg}; thus, let $0,T(0),T^2(0),\ldots,T^n(0)$ be
  the elements of the attractor ${\cal A}$, where we have $T^{n+1}(0)=0$.
  The result follows immediately, using the same argument as in the proof of
  Lemma \ref{ThmChar}.
\qed \end{pf*}

\paragraph*{Examples.}
We give some examples of zero expansions, which we write starting from the
least significant digit.
\begin{romanlist}
  \parskip=\smallskipamount
  \smallskip
  \item
    If $0\in {\cal D}$, the zero expansion is simply $(0)$.
  \item
    Take a base $b\in\Z$, with $|b|\ge 2$, and take digits
    $\{1,2,\ldots,|b|\}$. Obviously, the zero digit here is $|b|$. In this
    case, we have a zero expansion if and only if $b<0$. Indeed, if $b<0$, the
    zero expansion is given by $(|b|,1)$, because
    \begin{displaymath}
      |b|\cdot b^0 + 1\cdot b^1 = 0.
    \end{displaymath}
    If $b>0$, we cannot have a zero expansion: we have $d(0)=b$, so $T(0)=
    (0-b)/b=-1$, but negative numbers cannot be represented by nonnegative
    digits on a positive base. Indeed, we have $T(-1)= \frac{-1-(b-1)}{b}=-1$,
    so the zero expansion would be the infinite sequence
    $(b,b-1,b-1,b-1,\ldots)$. This implies immediately that 
    \begin{displaymath}
      (\Z,\,b,\, \{1,\,\ldots,\,|b|\}\,)
    \end{displaymath}
    for $b>0$ cannot be a number system.
  \item
    If $b\ge 2$, and we take the digits $\{-1,1,2,\ldots,b-2,b\}$, then we have
    the zero expansion $(b,-1)$. Note that this digit set gives a number system
    for any $b$, by Theorem \ref{ThmSmall}.
  \item
    We will show in Theorem \ref{ThmZeroExpansionLength} that the length of the
    zero expansion increases with the size of the zero digit. As an example of
    this behaviour, let $b=-2$, choose an integer $i\ge 0$, and let ${\cal D} =
    \{ 1, 3^i+1\}$; by Theorem \ref{ThmMinus2} below, this always gives a
    number system. The zero digit here is the even number $3^i+1$; it follows
    from Lemma \ref{Lemma3} that the zero expansion has length $3^i$.
\end{romanlist}

For a general digit set, the length of the zero expansion becomes an important
parameter in many kinds of number system constructions. For example, if we want
to pad an expansion to obtain some exact length $\ell$, we must know that the
length to be padded is divisible by the length of the zero expansion. This
problem will occur in the proof that there are infinitely many digit sets not
containing zero, for any base $b\in\Z$ (Theorems \ref{ThmAvoid} and
\ref{ThmAvoidNeg} below).

The last result in this subsection shows that in general, the length of the
zero expansion grows to infinity with the size of the zero digit.

\begin{thm} \label{ThmZeroExpansionLength}
  Let $(\Z,b,{\cal D})$ be a pre-number system with $0\in{\cal D}$. Then for
  each $\ell\ge 1$, there are only finitely many $d\in b\Z$ such that the
  pre-number system $(\Z,b,{\cal D}\backslash \{0\} \cup \{ d \} )$ has a
  zero expansion of length $\ell$.
\end{thm}

\begin{pf*}{Proof.}
  Let $\ell\ge 1$, let $d\in b\Z$, and let $(d_0,d_1,\ldots,d_{\ell-1})$ be
  the zero expansion of $(\Z,b,{\cal D}\backslash \{0\} \cup \{d\})$, as
  defined by \eqref{EqZeroCycle}. Let $I$ be the set of those $i$ in 
  $\{0,\ldots,\ell-1\}$ for which $d_i=d$; note that $0\in I$, because $d_0=d$.
  From $\sum_{i=0}^{\ell-1} d_i b^i=0$, we then obtain
  \begin{equation} \label{EqDFinite}
    \left( \sum_{i\in I} b^i \right)d = -\sum_{i\not\in I} b^i d_i.
  \end{equation}
  The element on the left is nonzero, because $|b|\ge 2$, and hence a sum of
  distinct powers of $b$ cannot be $0$.

  Now we finish the proof of the Theorem. The right hand side of 
  \eqref{EqDFinite} clearly takes at most $(|b|-1)^{\ell-1}$ distinct values.
  To each of these values corresponds at most one value for $d$.
  This completes the proof.
\qed \end{pf*}

\subsection{The first digit sets}

Note that the base $b=2$, although it can be used to define pre-number systems,
must be excluded. In fact, $b-1=1$ in this case, and Corollary
\ref{CorOneCycle} then tells us that there exists no digit set $\{d_0,d_1\}$ in
$\Z$ such that $(\Z,2,\{d_0,d_1\})$ is a number system. For example, the
well-known binary digits $\{0,1\}$ can only represent nonnegative integers on
base $2$.

The restriction to just $2$ digits is important here: for example, one can show
that every integer has a unique Non-Adjacent Form (NAF) expansion on base $2$
with the digits $\{0,1,-1\}$ (see \cite{Bosma2001}). In formulae: every $a\in
\Z$ can be written uniquely in the form
\begin{displaymath}
  a=\sum_{i=0}^\ell d_i2^i, \qquad d_i\in\{0,\pm 1\}, d_id_{i+1}=0.
\end{displaymath}
In this paper, however, we only consider irredundant digit sets, hence only
digit sets of cardinality $|b|$ if the base is $b$. Our first result here,
which is new as far as digit sets without zero are concerned, is as follows.

\begin{thm} \label{ThmSmall}
  Let $b\in\Z$, with $|b|\ge 3$. Let ${\cal D}$ be a complete residue system 
  modulo $b$, such that
  \begin{romanlist}
    \item $|d| \le |b|$ for all $d\in{\cal D}$;
    \item either $1\in {\cal D}$ or $-1\in{\cal D}$;
    \item neither $b-1\in {\cal D}$ nor $-b+1\in{\cal D}$.
  \end{romanlist}
  Then $(\Z,b,{\cal D})$ is a number system.
\end{thm}

\begin{pf*}{Proof.}
  Define $T : \Z\rightarrow\Z$ as in Definition \ref{DefT}; by Theorem
  \ref{ThmChar}, we must prove that for all $a\in\Z$, there exists $n\ge 1$
  such that $T^n(a)=0$. By the Lemma, it is enough to do this for all $a$ with
  $|a|\le 2$, as $|k|\le |b|$ and $|K|\le |b|$ in our case.

  For any $b$, if $|a|=1$, we easily verify that either $Ta=0$ or $T^2a=0$,
  using the second and third assumptions. If $a=2$, then either $2$ or $-|b|+2$
  is a digit, so that $T(2)\in\{0,1,-1\}$, and the same holds for $a=-2$.
  
  We see that for all nonzero $a\in\Z$, we have $T^na=0$ for some $n$.
  This immediately also shows the existence of a zero cycle, because if
  $a=T(0)$, there exists $n\ge 0$ such that $T^na=0$, so that
  $T^{n+1}(0)=0$.
\qed \end{pf*}

\paragraph*{Remarks.}
Note that the proof actually allows to relax condition (i) to $|d|\le 
2|b|-2$.

The above result does not hold as stated for base $-2$. Base $-2$ is actually
a quite special case, which will be worked out completely in Section
\ref{SecMinus2}.

The assumptions about the presence of $\pm 1$ in ${\cal D}$ are necessary. The
only representatives of $\pm 1$ that are allowed are $\pm 1$ themselves, $b\pm
1$, and $-b\pm 1$. If $b-1$ or $-b+1$ are digits, then we get a nonzero
$1$-cycle by Lemma \ref{Lem1Cycle}. If both $b+1$ and $-b-1$ are in ${\cal D}$,
we see 
\begin{displaymath}
  T(1)=\frac{1-(b+1)}{b} = -1, \quad 
  T(-1)=\frac{-1-(-b-1)}{b} = 1,
\end{displaymath}
which also gives a non-zero cycle.

\paragraph*{Example.}
A nice example of a digit set without zero that always works, is given by
the \emph{odd digit set}.

\begin{defn} \label{DefOdd}
For an \emph{odd} $b\in\Z$, define the set of \emph{odd digits modulo $b$} as
\begin{displaymath}
  {\cal D}_{b,\rm odd} = \begin{cases}
    \{\,-b+2,\,-b+4,\,\ldots,\,-1,\,1,\,\ldots,\,b-2,\,b\,\} & \text{ if } b>0;
    \\ 
    \{\, b,\,b+2,\,\ldots,\,-1,\,1,\,\ldots,\,-b-2\,\} & \text{ if } b<0.
    \end{cases}
\end{displaymath}
\end{defn}

\begin{cor} \label{CorOdd}
  Let $b\in\Z$ be odd, with $|b|\ge 3$. Then $(\Z,b,{\cal D}_{b,\rm odd})$ is
  a number system.
\end{cor}

\begin{pf*}{Proof.}
  The only thing to show, before we can apply the Theorem, is that ${\cal
  D}_{b,\rm odd}$ contains a complete system of representatives modulo $b$.
  Now if $d\in {\cal D}_{b,\rm odd}$ is negative, then $d+b$ is even and
  between $0$ and $b-1$. Thus the classes of $\{0,\ldots,b-1\}$ or
  $\{1,\ldots,b\}$ are all represented in ${\cal D}_{b,\rm odd}$.
\qed \end{pf*}

\subsection{Translation of digit sets}

In the quest for classification of all valid digit sets, now that we know that
having $0$ as a digit is not essential, we might think that one valid digit set
could give rise to infinitely many digit sets by simple translation. Below, we
show (Theorem \ref{ThmTranslate}) that translation of the digit set over a
fixed integer will destroy the number system property if the integer is too
large. In fact, we prove that when $0\in {\cal D}$, the same holds if we
translate all nonzero digits, while leaving $0$ in place. 

We begin with a basic observation.

\begin{lem} \label{LemDigAttr}
  Let $(\Z,b,{\cal D})$ be a pre-number system, with attractor ${\cal A}$. If
  $(\Z,b,{\cal D})$ is a number system, then ${\cal A}$ contains at least one
  element of ${\cal D}$.
\end{lem}

\begin{pf*}{Proof.}
  Consider the zero cycle 
  \begin{displaymath}
    0\rightarrow a_1\rightarrow \ldots \rightarrow a_\ell \rightarrow 0,
  \end{displaymath}
  where $a_1,\ldots,a_\ell\in{\cal A}$. If $a_\ell\rightarrow 0$, that means
  that the digit representing the coset of $a_\ell$ is equal to $a_\ell$, in
  other words, that $a_\ell\in{\cal D}$.
\qed \end{pf*}

The elements of the attractor may be thought of as ``small'', at least when
compared the the largest digit; therefore, the previous lemma tells us that at
least one digit is ``small''. However, we want to strengthen the claim of the
lemma to say that at least one \emph{nonzero} digit must be small. Note that
when $0\in{\cal D}$, the number system property is equivalent to ${\cal
A}=\{0\}$, so this nonzero digit cannot be an element of the attractor. The
next result, which generalises Theorem 4 from \cite{KoPe1983}, shows that next
to the attractor also the set $\{ a \in \Z \mid |a| \le L \}$ from Lemma
\ref{LemBd} has some importance.

\begin{lem} \label{LemKovacs}
  Let $(\Z,b,{\cal D})$ be a number system, and let $K$ and $L$ be as in Lemma
  \ref{LemBd}. If $K$ is large enough, then there is at least one $d\in{\cal
  D}$, with $d\ne 0$, such that 
  \begin{displaymath}
    |d| \le L.
  \end{displaymath}
\end{lem}

\begin{pf*}{Proof.}
  Let $a_0\in \Z$ have $a_0\ne 0$ and $|a_0|\le L$; we may assume that $K$ is
  so large that $L\ge 1$. By our assumption, $a_0$ has a finite expansion on
  the base $b$ with digits in ${\cal D}$. Thus, there exist a minimal $\ell$
  and $a_i\in \Z$ with 
  \begin{displaymath}
    a_0\rightarrow T(a_0)=a_1 \rightarrow \cdots \rightarrow a_\ell
    \rightarrow 0.
  \end{displaymath}
  By Lemma \ref{LemBd}, we know that $|a_i|\le L$ for all $i$. On the other
  hand, $a_\ell$ must be a digit, and by the minimality of $\ell$ we know that
  $a_\ell\ne 0$.
\qed \end{pf*}

\begin{thm} \label{ThmTranslate}
  Let $(\Z,b,{\cal D})$ be a pre-number system with $|b|\ge 3$, and for $t\in
  \Z$, define ${\cal D}_t$ as $\{ d+t \mid d\in{\cal D} \}$. Then there are
  only finitely many $t\in \Z$ such that $(\Z,b,{\cal D}_t)$ is a number
  system. \\
  If $0\in {\cal D}$, then the same statement holds for ${\cal \tilde{D}}_t=
    \{d+t \mid d\in {\cal D}, d\ne 0\} \cup \{0\}$.
\end{thm}

\begin{pf*}{Proof.}
  Let $K_t=\max_{d\in{\cal D}_t} |d|$; by Lemma \ref{LemBd}, we see that
  \begin{equation} \label{EqDigBd}
    |a| \le K_t /(|b|-1)
  \end{equation} 
  for all $a$ in the attractor ${\cal A}_t$ of $(\Z,b,{\cal D}_t)$.
  In particular, by Lemma \ref{LemDigAttr}, this inequality holds for at least
  one digit in ${\cal D}_t$; note that $1/(|b|-1)<1$ by our assumptions.
  But as $|t|\rightarrow\infty$, clearly $\frac{|d|}{K_t}\rightarrow 1$ for
  all $d\in{\cal D}_t$, so that \eqref{EqDigBd} is violated for all $d\in{\cal
  D}_t$ when $|t|$ is sufficiently large. This is a contradiction, and the
  first claim is proved.

  For the second claim, we use Lemma \ref{LemKovacs} to show that, when $t$ is
  large enough, we must have $|d| \le K_t / (|b|-1)$ for some nonzero
  $d\in{\cal \tilde{D}}_t$. The rest of the argument is the same.
\qed \end{pf*}

\paragraph*{Remark.}
The argument of the proof makes essential use of the inequality
$1/(|b|-1)<1$, and therefore the proof breaks down when $|b|=2$.  In fact, we
will obtain the assertions of the Theorem for the case $|b|=2$ below, using a
specialised argument.

\section{Infinitely many digit sets}

Having established the existence of good digit sets with and without zero for
any integer base $b$ with $|b|\ge 2$ (except $b=2$) in Theorem \ref{ThmSmall},
we will now proceed to show that every base $b$ with $|b|\ge 4$ has
\emph{infinitely many} good digit sets with and without zero --- see
Corollaries \ref{CorInf} and \ref{CorInfNeg}. This was already shown for digit
sets with zero by B. Kov\'acs and A. Peth\H o \cite[Section 4]{KoPe1983} for
negative $b$, and by D. Matula \cite{Matula} for any integer $b$ (both taking
$|b|\ge 3$). We will generalise their methods to our case.

Lemmas \ref{LemDigAttr} and \ref{LemKovacs} tell us that at least one nonzero
digit in the set must be small. The approach of Kov\'acs and Peth\H o is to
start from the standard digits $\{0,1,\ldots, |b|-1\}$ and replace just one
digit by a much larger number. We will adapt their proof to start from any
good digit set such that $|d|\le |b|$ for all digits $d$, and use this to show
that for any integer base $b$ with $|b|\ge 4$ there exist infinitely many digit
sets ${\cal D}$, both with and without zero, such that $(\Z,b,{\cal D})$ is a
number system. The case $|b|=3$ unfortunately remains open, as our methods do
not work for it. For the special case $b=-2$ we will characterise \emph{all}
valid digit sets later (see Theorem \ref{ThmMinus2} below).

\begin{defn}
Let $(\Z,b,{\cal D})$ be a pre-number system. If $a=\sum_{i=0}^\ell d_i b^i$
for some digits $d_i\in{\cal D}$, we say that $a$ has \emph{length $\ell+1$},
and write $L(a)=\ell+1$. Assume $a\ne 0$. If the expansion for $a$ is minimal,
and therefore unique, we call $d_\ell$ the \emph{most significant digit} of
$a$, and write $\MSD(a)$.
\end{defn}

Note that we have $\MSD(a)\ne 0$ by the minimality assumption.

Besides the functions $L(a)$ and $\MSD(a)$, we will use the following notation.
We let $(\Z,b,{\cal D})$ be a number system, such that $|d|\le |b|$ for all
$d\in{\cal D}$. We fix some $u\in\Z$ with $1\le |u|\le |b|-1$, some integer
$k\ge 1$, and one digit $d\in{\cal D}$, which is not the zero digit. Then, let
$\tilde{d} = d-ub^k$, and $\tilde{\cal D}= {\cal D} \setminus \{d\} \cup \{
\tilde{d} \}$. We write ${\cal A}$ and $\tilde{\cal A}$ for the attractors of
$(\Z,b,{\cal D})$ and $(\Z,b,\tilde{\cal D})$, respectively. 

\paragraph*{The case where $b>0$.}
We want to derive conditions on $u$ and $d$ that allow us to conclude that
$(\Z,b,\tilde{\cal D})$ is a number system for infinitely many values of $k$.
We start with a sharp lower bound on numbers with a given expansion length.
Recall that we assume $|d|\le b$ for all $d\in{\cal D}$.

\begin{lem} \label{LemLength}
  Assume $b\ge 3$, and let $a=\sum_{i=0}^\ell d_ib^i$ be a minimal expansion,
  with digits in ${\cal D}$. Then $a$ and $d_\ell$ have the same sign, and:
  \begin{romanlist}
    \item if $0\in{\cal D}$, then 
          $|a|\ge \frac{b^\ell+b-2}{b-1}$;
    \item if $0\notin{\cal D}$, then 
          $|a|\ge \frac{b^\ell-2b^{\ell-1}+b}{b-1}$.
  \end{romanlist}
\end{lem}

\begin{pf*}{Proof.}
  As $b>0$, by Theorem \ref{ThmSmall} and the remarks following it, both $1$
  and $-1$ are in ${\cal D}$, while neither $b-1$ nor $-b+1$ are in ${\cal D}$.
  Thus, we have $|d|\le b-2$ whenever $d\not\equiv 0\pmod{b}$.

  Suppose that $0\in{\cal D}$. Then we know that $|d_i|\le b-2$ for all
  $i$, and therefore
  \begin{displaymath}
    \left| {\textstyle \sum_{i=0}^{\ell-1} d_ib^i } \right| \le
    (b-2)\tfrac{b^\ell-1}{b-1} < b^\ell.
  \end{displaymath}
  Furthermore, we have $d_\ell\ne 0$ by the minimality assumption. It follows
  that
  \begin{displaymath}
    |a| \ge b^\ell - (b-2) \tfrac{b^\ell-1}{b-1} = \tfrac{b^\ell+b-2}{b-1}.
  \end{displaymath}

  If $0\notin{\cal D}$, minimality means that the expansion does not start
  with the zero expansion $(b,-1)$ or $(-b,1)$ (depending on whether $b$ or
  $-b$ is in ${\cal D}$). Thus, either $|d_\ell|>1$ or $|d_\ell|=1$ and 
  $|d_{\ell-1}|\le b-2$. In the first case, $|d_i|\le b$ for $0\le
  i\le \ell-1$, so
  \begin{displaymath}
    |a| \ge 2b^\ell - b \tfrac{b^\ell-1}{b-1} = \tfrac{b^{\ell+1}-1}{b-1}.
  \end{displaymath}
  In the second, we have
  \begin{displaymath}
    \left| {\textstyle \sum_{i=0}^{l-1} d_ib^i } \right| \le
    (b-2) b^{\ell-1} + b\left(b^{\ell-2}+\ldots+b+1\right) =
    b^\ell - 2b^{\ell-1} + \tfrac{ b^\ell-b }{b-1} <
    b^\ell,
  \end{displaymath}
  so that
  \begin{displaymath}
    |a| \ge b^\ell - \left( b^\ell - 2b^{\ell-1} + \tfrac{ b^\ell-b }{b-1}
    \right) = \tfrac{ b^\ell - 2b^{\ell-1} + b }{b-1}.\qed
  \end{displaymath}
\end{pf*}

\begin{lem} \label{LemU}
  Assume $b\ge 3$; if $0\not\in {\cal D}$, also assume $|u|\le b-2$. If $a$ is
  in $\tilde{\cal A}$, then $L(a)\le k+1$, and $L(a)=k+1$ implies
  $|\MSD(a)|= 1$.
\end{lem}

\begin{pf*}{Proof.}
  Let $a\in \tilde{\cal A}$. We may assume that $\tilde{d} = d-ub^k$ has
  maximum absolute value in $\tilde{\cal D}$, since otherwise $|\tilde{d}|\le
  b$ and we can apply Theorem \ref{ThmSmall} to decide if $\tilde{\cal D}$ is
  a valid digit set. Thus by Lemma \ref{LemBd}, we have $|a| \le
  \frac{|u|b^k+|d|}{b-1} \le \frac{|u|}{b-1}b^k + 1$.

  If $0\in{\cal D}$, this bound is simply $|a|\le b^k+1 = \frac{b^{k+1}-b^k+b
  -1}{b-1}$. Now assume also that $L(a)\ge k+2$; then by Lemma \ref{LemLength},
  we have $|a|\ge \frac{b^{k+1}+b-2}{b-1}$. This is a contradiction.

  If $0\not\in{\cal D}$, we assumed $|u|\le b-2$, so
  $|a|\le \frac{ b^{k+1} - 2b^k + b-1}{b-1}$. Assume that $L(a)\ge k+2$; then by
  Lemma \ref{LemLength}, we have $|a|\ge \frac{b^{k+1}-2b^k+b}{b-1}$, which is
  impossible.

  Now assume $L(a)=k+1$, and $|\MSD(a)|>1$. Then the lower bounds for $|a|$
  given by Lemma \ref{LemLength} are $b^k + \frac{b^k+b-2}{b-1}$ and
  $b^k + \frac{b^k-2b^{k-1}+b}{b-1}$, respectively, and these are still in
  contradiction with $|a|\le \frac{|u|}{b-1}b^k+1$.
\qed \end{pf*}

\begin{lem} \label{LemMSD}
  Assume $b\ge 3$. Let $a\in\Z$ have $|a|\le b-1$; then $L(a)\le 2$, and if
  $L(a)=2$, then $|\MSD(a)|= 1$.
\end{lem}

\begin{pf*}{Proof.}
  This follows directly from Lemma \ref{LemLength}: if we assume $L(a)=3$,
  we find $|a|\ge b$, a contradiction, and the same happens if we assume
  $L(a)=2$ and $|\MSD(a)|\ge 2$.
\qed \end{pf*}

\begin{defn} \label{DefDk}
Assume $b\ge 3$. Recall our fixed digit $d\in{\cal D}$.
For an integer $k\ge 0$, define
\begin{align*}
  D_k &= \{ (d_0,d_1,\ldots,d_k) : d_i \in {\cal D} \text{ for }
  0\le i\le k-1, \, d_k\in \{-1,0,1\} \}, \\
  \tilde{D}_k &= \{ (d_0,d_1,\ldots,d_k)\in D_k : d_i\ne d \text{ for }
  0\le i\le k \}.
\end{align*}
\end{defn}

Clearly, $D_k$ contains all digit expansions with digits in ${\cal D}$ and
length padded to exactly $k+1$, such that the most significant digit is at
most $1$ in absolute value. If $0\not\in{\cal D}$, we still allow $d_k=0$,
because otherwise it is not always possible to pad exactly to the required 
length. The subset $\tilde{D}_k$ consists of all elements of $D_k$ that have
no components equal to $d$.

Next, we define the function $\Phi_k : D_k\rightarrow D_k$ as follows. Let
${\bf d}=(d_0,\ldots,d_k)\in D_k$. If $d_0=d$, our fixed digit, then
\begin{subequations}
\begin{align}
  \Phi_k({\bf d}) &= (d_1,\ldots,d_{k-1},d_0',d_1') \label{EqPhi1} \\
  \intertext{where $d_0'$ and $d_1'$ in ${\cal D}$ are such that 
    $d_0'+d_1'b=d_k+u$. This is possible by Lemma \ref{LemMSD}. If $d_0\ne d$,
    then} 
  \Phi_k({\bf d}) &= \begin{cases}
    (d_1,\ldots,d_k,0) & \text{ if } d_k\ne 0 \text{ or } 0\in {\cal D} \\
    (d_1,\ldots,d_{k-1},d_0',d_1') & \text{ otherwise,}
      \label{EqPhi2}
    \end{cases}
\end{align}
\end{subequations}
where $d_0'$ and $d_1'$ in ${\cal D}$ satisfy $d_0'+d_1'b=0$.

\begin{lem} \label{LemCriterion}
  Assume $b\ge 3$; if $0\not\in{\cal D}$, also assume $|u|\le b-2$. Then
  $\Phi_k$ is well defined. Furthermore, if for each ${\bf d}\in D_k$ there
  exists an $n\ge 0$ such that $\Phi_k^n({\bf d}) \in \tilde{D}_k$, then
  $(\Z,b,\tilde{\cal D})$ is a number system.
\end{lem}

\begin{pf*}{Proof.}
  We extend an argument that was already used in \cite{KoPe1983,Matula}. It
  runs as follows. In order to prove that $(\Z,b,\tilde{\cal D})$ is a number
  system, it suffices that every element in
  the attractor $\tilde{\cal A}$ has a finite expansion with digits in
  $\tilde{\cal D}$.  Let $a\in\tilde{\cal A}$ and let $a=\sum_{i=0}^k d_i b^i$
  be its expansion with digits in ${\cal D}$, padded to length $k+1$; if
  necessary, we may take $d_k=0$, even if $0\not\in {\cal D}$. It follows that
  ${\bf d}=(d_0,\ldots,d_k)$ is in $D_k$.

  There are two cases. If $d_0\ne d$, then $a$ has a finite expansion with
  digits in $\tilde{\cal D}$ if and only if $(a-d_0)/b$ has such an expansion.
  If $d_0=d$, we replace $d_0$ by $d-ub^k$; to make up, we also replace $d_k$
  by $d_k+u$. Then $a$ has an expansion of the desired form if and only if
  $(a-(d-ub^k))/b$ does.

  We claim that if ${\bf d}\in D_k$ is an expansion of $a$, then 
  $\Phi_k({\bf d})$ is an expansion of $(a-d_0)/b$ and $(a-(d-ub^k))/b$,
  in the respective cases. Clearly, if this claim holds, then the lemma follows
  by induction, because the expansions in $\tilde{D}_k$ are of the desired
  form.

  We prove the claim, using the same two cases. Let ${\bf d}\in D_k$ be an
  expansion of $a$, with $d_0\ne d$; then we get an expansion of $(a-d_0)/b$ by
  deleting $d_0$ and shifting the other digits down. To have an expansion of
  length $k+1$ again, we can add a digit $0$ if $0\in {\cal D}$. If
  $0\not\in{\cal D}$, we must be careful. If $d_k=0$ already, by adding $0$ we
  would get two zeros in succession, and this is not allowed by the definition
  of $D_k$. Instead, we also delete $d_k$, and add the zero expansion of
  $(\Z,b,{\cal D})$, which is either $(b,-1)$ or $(-b,1)$. If $d_k\ne 0$,
  however, we cannot do this, and we add a $0$. This corresponds to the
  definition \eqref{EqPhi2} of $\Phi_k({\bf d})$ in this case.

  If $d_0=d$, as already said, we replaced $d_0$ by $d-ub^k$, and $d_k$ by
  $d_k+u$. Now consider $(a-(d-ub^k))/b$. As before, we delete $d-ub^k$ and
  shift the other digits down. Of course, $d_k+u$ need not be a digit. However,
  because $d_k\in\{-1,0,1\}$ and $|u|\le b-1$ or $b-2$, according as $0\in{\cal
  D}$ or not, $d_k+u$ can be written $d_0'+d_1'b$ with $d_0'\in{\cal D}$ and
  $d_1'\in\{-1,0,1\}$, by Lemma \ref{LemMSD}. Therefore, we replace $d_k+u$
  by this expansion of length $1$ or $2$, adding a $0$ if necessary. This gives
  us an expansion of length $k+1$ that satisfies the definition of $D_k$, and
  corresponds to the definition of $\Phi_k({\bf d})$ in \eqref{EqPhi1}. The
  claim is proved.
\qed \end{pf*}

The next result generalises Theorem 5 in \cite{KoPe1983}.

\begin{thm} \label{ThmAvoid}
  Let $(\Z,b,{\cal D})$ be a number system, where $b\ge 3$ and where
  $|d|\le b$ for all $d\in{\cal D}$. Fix some $d\in{\cal D}$ and some integer
  $u$ with $|u|\le b-1$; if $0\not\in{\cal D}$, assume $|u|\le b-2$. Let
  ${\cal B}$ be the set of digits in ${\cal D}$ that occur in the expansions of
  $0$, $u+1$, $u$, and $u-1$. If $d\not\in{\cal B}$, then we may replace $d$
  in ${\cal D}$ by $\tilde{d}=d-ub^k$, for any $k\ge 1$, without affecting the
  number system property.
\end{thm}

\begin{pf*}{Proof.}
  Let ${\bf d}\in D_k$, as defined above; by Lemma \ref{LemCriterion}, it is
  enough to show that $\Phi_k({\bf d}) \in \tilde{D}_k$ for $n$ large enough.
  Now whatever the components of ${\bf d}$ are, they are gradually replaced by
  the components introduced at the end by the repeated application of $\Phi_k$.
  These new components are digits that occur in the expansion of $0$, of $1+u$,
  of $u$, and of $-1+u$. Thus if $d$ is distinct from all these digits, then
  for $n$ large enough, $\Phi_k({\bf d})$ will have no components equal to $d$,
  as desired.
\qed \end{pf*}

\paragraph*{Remarks.}
  The least significant digits of $0$, $u$, $u+1$, and $u-1$, and the possible
  most significant digits $1$ and $-1$, together make up the set ${\cal B}$.
  Therefore, ${\cal B}$ has at most $6$ elements.

  It follows from the proof that the zero digit, being the least significant
  digit of $0$, is always one of the bad digits, and in fact the conclusion of
  the Theorem is often false if $d$ is congruent to $0$ modulo $b$. For
  example, although ${\cal D}=\{-5,\,1,\,2,\,3,\,-1\}$ gives a number system
  with base $b=5$, the sets $\{-5+5^k,\,1,\,2,\,3,\,-1\}$ for $k\ge 2$ give a
  cycle $(5-5^k)/4\rightarrow (5-5^k)/4$, and the attractors of
  $\{-5-5^k,\,1,\,2,\,3,\,-1\}$ for $k\ge 2$ do not contain $0$.

\paragraph*{Examples.}
  Let us apply Theorem \ref{ThmAvoid} to some of the starting digit sets that
  we found in the previous section.

  First, let us note that Theorem \ref{ThmAvoid} cannot be applied if $b=3$.
  Indeed, because $u$, $u+1$, and $u-1$ are incongruent modulo $b$, we see
  that ${\cal B}$ must contain at least $3$ elements. If now $b=3$, we have
  no choices left for $d$. 

  In fact, we have been unable to find any infinite sequence of valid nonzero
  digit sets for $b=3$. However, the set $\{0,1,2-3^k\}$ was found to be valid
  for all $k\ge 1$ by Matula \cite[Theorem 8]{Matula}. He used a refinement of
  our argument for the case where ${\cal D}$ has only nonnegative digits, which
  allows him to start from the digit set $\{0,1,2\}$. Of course, with this
  digit set only nonnegative integers can be represented, but using Theorem
  \ref{ThmZBd} one can prove that the attractor $\tilde{\cal A}$ contains only
  nonnegative elements if we choose $u$ positive. Unfortunately, this argument
  does not work in the case the starting digit set contains $b$ instead of $0$.

  Due to these technical problems with $b=3$, we assume $b\ge 4$ in what
  follows.

  Consider ${\cal D}=\{-1,\,0,\,1,\,\ldots,\,b-2\}$; this is a valid digit set
  by Theorem \ref{ThmSmall}. Taking $u=1$, we find the expansions $0=(0)$,
  $u=(1)$, $u+1=(2)$, and $u-1=(0)$. It follows that ${\cal B}=\{0,1,2\}$,
  so we can take $d=-1$ or $d=3,\,4,\,\ldots,\,b-2$ and replace it by
  $\tilde{d}=d-b^k$ for any $k\ge 1$. If we take $u=-1$, the expansion for
  $u-1$ becomes $-2=(b-2,\,-1)$, and we obtain ${\cal B}=\{-1,\,0,\,b-2\}$. 
 
  For an example without the digit $0$, consider 
  ${\cal D}=\{-1,\,1,\,2,\,\ldots,\,b-2,\,b\}$. Again by Theorem
  \ref{ThmSmall}, this digit set is valid. Taking $u=1$, we expand
  $0=u-1=(b,\,-1)$, $u=(1)$, and $u+1=(2)$, so that ${\cal B}
  =\{-1,\,1,\,2,\,b\}$. For $u=-1$, we get ${\cal B}=\{-1,\,b-2,\,b\}$.
  
  We thus obtain the following basic result.

  \begin{cor} \label{CorInf}
    For each integer base $b\ge 4$ there exist infinitely many valid digit sets
    ${\cal D}$ containing $0$, and infinitely many valid digit sets without
    $0$.
  \end{cor}

  \begin{pf*}{Proof.}
    For any $k\ge 1$, one can take $\{0,\,1,\,\ldots,\,b-2\}\,\cup\,\{
    -1-b^k\}$ and $\{-1,\,2,\,3,\,\ldots,\,b-2,\,b\}\,\cup\,\{ 1+b^k\}$,
    respectively.
  \qed \end{pf*}

  As another example, let $b$ be odd, and consider the odd digit set ${\cal
  D}_{\rm odd}$ (Definition \ref{DefOdd}). Let us choose $u=1$; we find the
  expansions $0=(b,\,-1)$, $u=(1)$, $u-1=(b,\,-1)$, and $u+1=(-b+2,\,1)$.
  Consequently, the bad set ${\cal B}$ is $\{1,\,-1,\,b,\,-b+2\}$. For $u=-1$,
  we have $u+1=(b-2,\,-1)$, and ${\cal B}=\{-1,\,b-2,\,b\}$.

\paragraph*{The case where $b<0$.} We now change to the case where the base $b$
is negative, still assuming that we start from a digit set ${\cal D}$ with all
digits at most equal to $|b|$ in absolute value. Obtaining upper bounds on the
expansion length is trickier here than before, because of the sign alternation
in powers of $b$ in consecutive terms of the expansion. The results are as 
follows. Note that we exclude $b=-2$; for this very special case, we refer to
Section \ref{SecMinus2} below.

\begin{lem} \label{LemLengthNeg}
  Assume $b\le -3$, and let $a=\sum_{i=0}^\ell d_ib^i$ be a minimal expansion, 
  with $d_i\in {\cal D}$. Let $L(0)$ be the length of the zero expansion with
  digits in ${\cal D}$. Then $a$ and $d_\ell b^{\ell}$ have the same sign,
  and, putting $B=|b|$, we have:
  \begin{romanlist}
    \item if $L(0)=1$ and $\ell\ge 0$, then $|a|\ge \begin{cases}
             1+\frac{B^\ell-1}{B^2-1} & \text{ if $\ell$ is even;} \\
	     1+\frac{B^\ell-B}{B^2-1} & \text{ if $\ell$ is odd.}
	     \end{cases}$
    \item if $L(0)=2$ and $\ell\ge 1$, then $|a|\ge \begin{cases}
             1+\frac{B^\ell-B^{\ell-1}+B-1}{B^2-1} 
	             & \text{ if $\ell$ is even;} \\
	     1+\frac{B^\ell-B^{\ell-1}-B+1}{B^2-1} & \text{ if $\ell$ is odd.}
	     \end{cases}$
    \item if $L(0)=3$ and $\ell\ge 2$, then $|a|\ge \begin{cases}
             1+\frac{B^\ell-2B^{\ell-2}+1}{B^2-1} 
	             & \text{ if $\ell$ is even;} \\
	     1+\frac{B^\ell-2B^{\ell-2}+B}{B^2-1} & \text{ if $\ell$ is odd.}
	     \end{cases}$
  \end{romanlist}
\end{lem}

\begin{pf*}{Proof.}
  We write $B=|b|$ throughout. As in the proof of Lemma \ref{LemLength}, we
  will show that $|\sum_{i=0}^{\ell-1} d_ib^i|$ is less than $B^\ell$ for
  minimal expansions. Thus, all claims will follow from the fact that
  \begin{displaymath}
    |a| \ge B^\ell - |{\textstyle \sum_{i=0}^{\ell-1} d_ib^i }|.
  \end{displaymath}
  Now minimising $|a|$ amounts to maximising the second term on the right.
  This can be done by maximising all $d_i$ with odd $i$, and minimising those
  with even $i$, or conversely. 

  First, assume $0\in{\cal D}$; this implies $|d|\le B-1$ for all $d\in{\cal
  D}$. Because ${\cal D}$ is a valid digit set, either $1$ or $-1$ is in 
  ${\cal D}$; let us assume the former. Thus the expansion with smallest
  absolute value is given by
  \begin{displaymath}
    (\ldots, \;\; b+2, \;\; -b-1, \;\;\; b+2, \;\; -b-1, \;\;\; 1).
  \end{displaymath}
  This is explained as follows: we take the most significant digit as small as
  possible, but cannot make it $0$ in a minimal expansion. Then we maximise
  the second digit, using something positive to get the sign right; we cannot
  get beyond $-b-1$. Then, we would like to take $b+1$ in the third digit,
  being maximally negative; but $b+1$ and $1$ cannot be in the same digit set.
  Thus, the third digit is $b+2$ or greater. We find that
  $\left|\,{\displaystyle \sum_{i=0}^{\ell-1} d_ib^i}\,\right|$ is bounded by 
  \begin{displaymath}
    (B-1)(B+B^3+\ldots+B^{\ell-1}) + (B-2)(1+B^2+\ldots+B^{\ell-2})
    = B^\ell-1 - \tfrac{B^\ell-1}{B^2-1}
  \end{displaymath}
  when $\ell$ is even, and by 
  \begin{displaymath}
    B^\ell-1-\tfrac{B^\ell-B} {B^2-1}
  \end{displaymath}
  when $\ell$ is odd.

  Next, assume we have a zero expansion of length $2$, which will be either
  $(-b,\; 1)$ or $(b,\; -1)$. Let us assume the former. Minimality now
  forbids to have $d_\ell=1$ and $d_{\ell-1}=-b$, so we may assume $d_\ell=1$
  and $d_{\ell-1}=-b-1$. Thus $1$, $-b$, and $-b-1$ are in ${\cal D}$, and we
  see that $b+1\not\in{\cal D}$. Therefore, the smallest expansion is given by
  \begin{displaymath}
    (\ldots, \;\; -b, \;\;\; b+2, \;\; -b, \;\;\; b+2, \;\; -b-1, \;\;\; 1).
  \end{displaymath}

  We find that $ \quad
    \left|\,{\displaystyle \sum_{i=0}^{\ell-1} d_ib^i}\,\right| \le B^\ell - 1 -
    \begin{cases}
      \tfrac{B^\ell-B^{\ell-1}+B-1}{B^2-1} &
        \text{ if $\ell$ is even;} \\[\smallskipamount]
      \tfrac{B^\ell-B^{\ell-1}-B+1}{B^2-1} &
        \text{ if $\ell$ is odd.}
    \end{cases} $

  Finally, assume the zero expansion is $(b,\; -b-1,\; 1)$ or $(-b,\; b+1,\;
  -1)$; let us say, the former. It follows that $-1$ and $b+1$ are not in
  ${\cal D}$, and the smallest expansion is given by
  \begin{displaymath}
    (\ldots, \;\; b, \;\; -b-1, \;\;\; b, \;\; -b-1, \;\;\; b+2, \;\; -b-1,
    \;\;\; 1).  
  \end{displaymath}
  We find that $ \quad
    \left|\,{\displaystyle \sum_{i=0}^{\ell-1} d_ib^i}\,\right| \le B^\ell - 1 -
    \begin{cases}
      \tfrac{B^\ell-2B^{\ell-2}+1}{B^2-1} &
        \text{ if $\ell$ is even;} \\[\smallskipamount]
      \tfrac{B^\ell-2B^{\ell-2}+B}{B^2-1} &
        \text{ if $\ell$ is odd. \qed}
    \end{cases} $
\end{pf*}

\begin{lem} \label{LemUNeg}
  Assume $b\le -3$, and write $B=|b|$; if $0\not\in {\cal D}$, assume $|u|\le
  B-2$. If $a$ is in $\tilde{\cal A}$, then $L(a)\le k+2$, and $L(a)>k$ implies
  that $|a-\sum_{i=0}^{k-1} d_ib^i|= B^k$.
\end{lem}

\begin{pf*}{Proof.}
  We write $B=|b|$ and let $a\in\tilde{\cal A}$. The method is the same as for
  Lemma \ref{LemU}, and we will leave the details to the reader. The fact that
  $a$ is in ${\cal A}$ leads to upper bounds on $|a|$, while lower bounds on
  $|a|$ are provided by Lemma \ref{LemLengthNeg}.

  The implication when $L(a)>k$ is proved as follows. If the implication is
  false, then the lower bounds from Lemma \ref{LemLengthNeg} for $\ell=k$ or
  $\ell=k+1$ can be increased by $B^k$, and this makes them larger than the
  upper bound for $|a|$.
\qed \end{pf*}

%


\begin{lem} \label{LemMSDNeg}
  Assume $b\le -3$, and write $B=|b|$. Let $a\in\Z$ with $|a|\le B-1$; then 
  $L(a)\le 3$, and if $L(a)>1$, then $|a-d_0|=B$.
\end{lem}

\begin{pf*}{Proof.}
  This follows directly from Lemma \ref{LemLengthNeg}: if we assume $L(a)=4$,
  we find $|a|\ge B$, a contradiction, and the same happens if we assume
  $L(a)=2$ or $L(a)=3$ and $|a-d_0|\ge 2B$.
\qed \end{pf*}

We now define a discrete dynamical system analogous to the one defined above;
see Definition \ref{DefDk}.
  
\begin{defn} \label{DefDkNeg}
Assume $b\le -3$. Recall our fixed digit $d\in{\cal D}$.
For an integer $k\ge 0$, define
\begin{align*}
  S &= \begin{aligned}[t]
          \{ (d_0,d_1) : {} & d_0,d_1\in{\cal D}\cup \{0\}, \; d_0+bd_1\in
             \{-1,0,1\}, \; \\
           & (d_0,d_1)\ne (0,0) \text{ if } L(0)=2 \};
        \end{aligned} \\
  E_k &= \{ (d_0,d_1,\ldots,d_{k+1}) : d_i \in {\cal D} \text{ for }
  0\le i\le k-1, \, (d_k,d_{k+1}) \in S \}; \\
  \tilde{E}_k &= \{ (d_0,d_1,\ldots,d_{k+1})\in D_k : d_i\ne d \text{ for }
  0\le i\le k+1 \}.
\end{align*}
\end{defn}

The set $E_k$ contains all expansions over ${\cal D}$ of length $k+1$ such that
the most significant part $d_k+bd_{k+1}$ has absolute value at most $1$. The
possible pairs $(d_k,d_{k+1})$ that satisfy this condition depend on ${\cal
D}$, and are collected in the set $S$. In order to get a length of exactly
$k+1$, we allow some digits to be $0$, even if $0$ is not in ${\cal D}$, just
as in the case $b>0$ (Definition \ref{DefDk}). Our definition implies that
$S$ has $3$ elements for every ${\cal D}$, namely the expansions of $-1$, $1$,
and $0$.

Next, we define the function $\Psi_k : E_k\rightarrow E_k$ as follows. Let
${\bf d}=(d_0,\ldots,d_{k+1})\in D_k$. If $d_0=d$, our fixed digit, then
\begin{subequations}
\begin{align}
  \Psi_k({\bf d}) &= (d_1,\ldots,d_{k-1},d_0',d_1',d_2') \label{EqPsi1} \\
  \intertext{where $d_0',d_1',d_2'$ in ${\cal D}$ are such that 
    $d_0'+d_1'b+d_2b^2=d_k+d_{k+1}b+u$. This is possible by Lemma
    \ref{LemMSDNeg}.  Suppose $d_0\ne d$. If $d_{k+1}\ne 0$ or $0\in {\cal D}$,
    then} 
  \Psi_k({\bf d}) &= (d_1,\ldots,d_{k+1},0). \label{EqPsi2}
  \intertext{If $d_{k+1}=0$ and $(d_0',d_1')$ is the zero expansion, then }
  \Psi_k({\bf d}) &= (d_1,\ldots,d_k,d_0',d_1'). \label{EqPsi3}
  \intertext{If $d_k=d_{k+1}=0$ and $(d_0',d_1',d_2')$ is the zero expansion,
  then}
  \Psi_k({\bf d}) &= (d_1,\ldots,d_{k-1},d_0',d_1',d_2'). \label{EqPsi4}
\end{align}
\end{subequations}

\begin{lem} \label{LemCriterionNeg}
  Assume $b\le -3$, and write $B=|b|$; if $0\not\in{\cal D}$, also assume 
  $|u|\le B-2$. Then $\Psi_k$ is well defined. Furthermore, if for each ${\bf
  d}\in E_k$ there exists an $n\ge 0$ such that $\Psi_k^n({\bf d}) \in
  \tilde{E}_k$, then $(\Z,b,\tilde{\cal D})$ is a number system.
\end{lem}

\begin{pf*}{Proof.}
  The fact that $\Psi_k$ is well defined, i.e., defines a map from $E_k$
  into $E_k$, follows directly from Lemma \ref{LemUNeg}. The rest of the
  argument is the same as for Lemma \ref{LemCriterion}. One uses
  Lemma \ref{LemMSDNeg} to show that $d_k+d_{k+1}b+u$ always has an 
  expansion of length at most $3$, so that $\Psi_k({\bf d})$ always ``fits''
  into the set $E_k$.
\qed \end{pf*}

\begin{thm} \label{ThmAvoidNeg}
  Let $(\Z,b,{\cal D})$ be a number system, where $b\le -3$, and where
  $|d|\le B$ for all $d\in{\cal D}$, with $B=|b|$. Fix some $d\in{\cal D}$ and
  some integer $u$ with $|u|\le B-1$; if $0\not\in{\cal D}$, assume $|u|\le
  B-2$. Let ${\cal B}$ be the set of digits in ${\cal D}$ that occur in the
  expansions of $0$, $u+1$, $u$, and $u-1$. If $d\not\in{\cal B}$, then we may
  replace $d$ in ${\cal D}$ by $\tilde{d}=d-ub^k$, for any $k\ge 1$, without
  affecting the number system property.
\end{thm}

\begin{pf*}{Proof.}
  Let ${\bf d}\in E_k$, as defined above; by Lemma \ref{LemCriterion}, it is
  enough to show that $\Psi_k({\bf d}) \in \tilde{E}_k$ for $n$ large enough.
  Now whatever the components of ${\bf d}$ are, they are gradually replaced by
  the components introduced at the end by the repeated application of $\Psi_k$.
  These new components are the digits that occur in the expansion of $0$, of
  $1+u$, of $u$, and of $-1+u$. Thus if $d$ is distinct from all these digits,
  then for $n$ large enough, $\Psi_k({\bf d})$ will have no components equal to
  $d$, as desired.
\qed \end{pf*}

\paragraph*{Remarks.}
  The same remarks as with Theorem \ref{ThmAvoid} apply here. The expansions
  of $0$, $u+1$, $u$, and $u-1$ among them have at most $4$ distinct least
  significant digits; the more significant digits $d_1$ and maybe $d_2$ are
  all taken from $\{1,\,-1,\,-b-1,\,b+1\}$. Therefore, $|{\cal B}|\le 8$.
  
  An example where the conclusion of the Theorem is false when $d\equiv
  0\pmod{b}$ is the following. Although ${\cal D}=\{-5,\,1,\,2,\,3,\,4\}$ gives a
  number system with base $b=-5$, the set $\{ -5-(-5)^k,\,1,\,2,\,3,\,4\}$ is not
  valid for $k\ge 2$: we have $-5-(-5)^k=-5(1-(-5)^{k-1})$, which is 
  divisible by $1-(-5)=6$ for $k\ge 1$, and thus gives a nonzero $1$-cycle 
  $\frac{-5-(-5)^k}6 \rightarrow \frac{-5-(-5)^k}6$ by Lemma \ref{Lem1Cycle} if
  $k\ge 2$.

\paragraph*{Examples.}
  Let $b<0$. For the reasons explained after Theorem \ref{ThmAvoid}, we cannot
  apply Theorem \ref{ThmAvoidNeg} when $b=-3$. Thus, assume $b\le -4$. We write
  $B=|b|$.

  Consider the classical digit set $\{0,\,1,\,\ldots,\,B-1\}$, and take $u=1$.
  It is clear that the bad set ${\cal B}$ is $\{0,1,2\}$, so we may replace $d$
  by $d-b^k$ for any $k\ge 1$, if $3\ge d\ge B-1$. Now take $u=-1$. We find
  $u=(B-1,\,1)$ and $u-1=(B-2,\,1)$, so that ${\cal B}=\{0,\,1,\,B-2,\,B-1\}$.
  Thus, any $d$ outside the latter set may be replaced by $d+b^k$, for any
  $k\ge 1$.

  Now as an example of a nonzero digit set, let ${\cal D}=\{1,2,\ldots,B\}$.
  We find $0=(B,\,1)$, as $B=-b$, and with $u=1$, we have $u=(1)$, $u-1=
  (B,\,1)$, and $u+1=(2)$. Thus ${\cal B}=\{1,\,2,\,B\}$. For $u=-1$, we find
  $u=(B-1,\,1)$, $u-1=(B-2,\,1)$, and $u+1=(B,\,1)$, so that ${\cal B}=
  \{1,\,B-2,\,B-1,\,B\}$.

  \begin{cor} \label{CorInfNeg}
    For each integer base $b\le -4$ there exist infinitely many valid digit
    sets ${\cal D}$ containing $0$, and infinitely many valid digit sets
    without $0$.
  \end{cor}

  \begin{pf*}{Proof.}
    For any $k\ge 1$, one can take $\{0,1,\ldots,B-2\}\,\cup\,\{ B-1-b^k\}$
    and $\{1,2,4,5,\ldots,B\}\,\cup\,\{ 3-b^k\}$, respectively.
  \qed \end{pf*}

  With the odd digits ${\cal D}_{\rm odd}$ (Definition \ref{DefOdd}), we
  have $0=(b,\,-1)$. For $u=1$, we get $u+1=(b+2,\,-1)$, so that ${\cal B}=
  \{-1,\,1,\,b+2,\,b\}$. For the more exotic $u=-3$, we get $u=(-3)$, 
  $u-1=(B-4,\,1)$, and $u+1=(B-2,\,1)$, so that ${\cal
  B}=\{1,\,-1,\,-3,\,B-4,\,B-2\}$.

  Finally, as an example of a digit set with a zero expansion of length $3$,
  let ${\cal D}=\{b,\,1,\,2,\,\ldots,\,B-1\}$ and $u=1$. This gives
  $0=(b,\,B-1,\,1)$, and ${\cal B}=\{ 1,\,2,\,B-1,\,b \}$.

~

It is an interesting question whether there also exist infinitely ``zero
digits'' complementing a given digit set. For example, for $b\le -2$, are there
infinitely many multiples $cb$ of $b$ such that $\{\,cb,1,2,\ldots,|b|-1\}$ is
a good digit set? As yet, we only have some partial answer to this question.
Namely, Theorem \ref{ThmZeroExpansionLength} shows that as $|c|\rightarrow
\infty$, with the other digits staying the same, also the length of the zero
cycle increases without bound. This contrasts with the infinite families that
we gave in this section, where the length of the zero cycle is the same 
throughout the family.

\section{Base $-2$} \label{SecMinus2}

The case where the base $b$ of the number system is $-2$ is special, as several
of the general results obtained above do not apply to this case. Examples are
Theorem \ref{ThmSmall} about smallest digit sets, Theorem \ref{ThmTranslate}
that says that only finitely many translates of a given digit set can yield
number systems, and the Theorems given in the last section that prove the
existence of infinitely many good digit sets.

However, in the case of the integers $\Z$, we have succeeded in determining 
\emph{all} possible digit sets for the base $b=-2$. It will follow from this
characterisation that there are infinitely many good digit sets for this base
and that unbounded translation only yields finitely many good such sets. A
remarkable feature of this case is that there exist no infinite families of
good digit sets obtained by translating one digit by a power of $-2$, as in the
last section; instead, one can shift by powers of $3$.

\begin{thm} \label{ThmMinus2}
  Let $d,D\in\Z$, with $d<D$. Then $(\Z,-2,\{d,D\})$ is a number system if and
  only if
  \begin{romanlist}
    \item one of $\{d,D\}$ is even and one is odd;
    \item neither $d$ nor $D$ is divisible by $3$, except that the even
          digit can be $0$;
    \item we have $2d\le D$ and $2D\ge d$;
    \item $D-d=3^i$ for some $i\ge 0$.
  \end{romanlist}
\end{thm}

~

As an example, the Theorem implies that a valid digit set for base $-2$ that
contains $0$ must be either $\{0,1\}$ or $\{0,-1\}$. On the other hand, it
follows easily that there are infinitely many valid digit sets without $0$,
for example the sets $\{1,3^i+1\}$ for $i\ge 0$ already discussed earlier.

\begin{figure*}[htp] \centering
  \epsfig{file=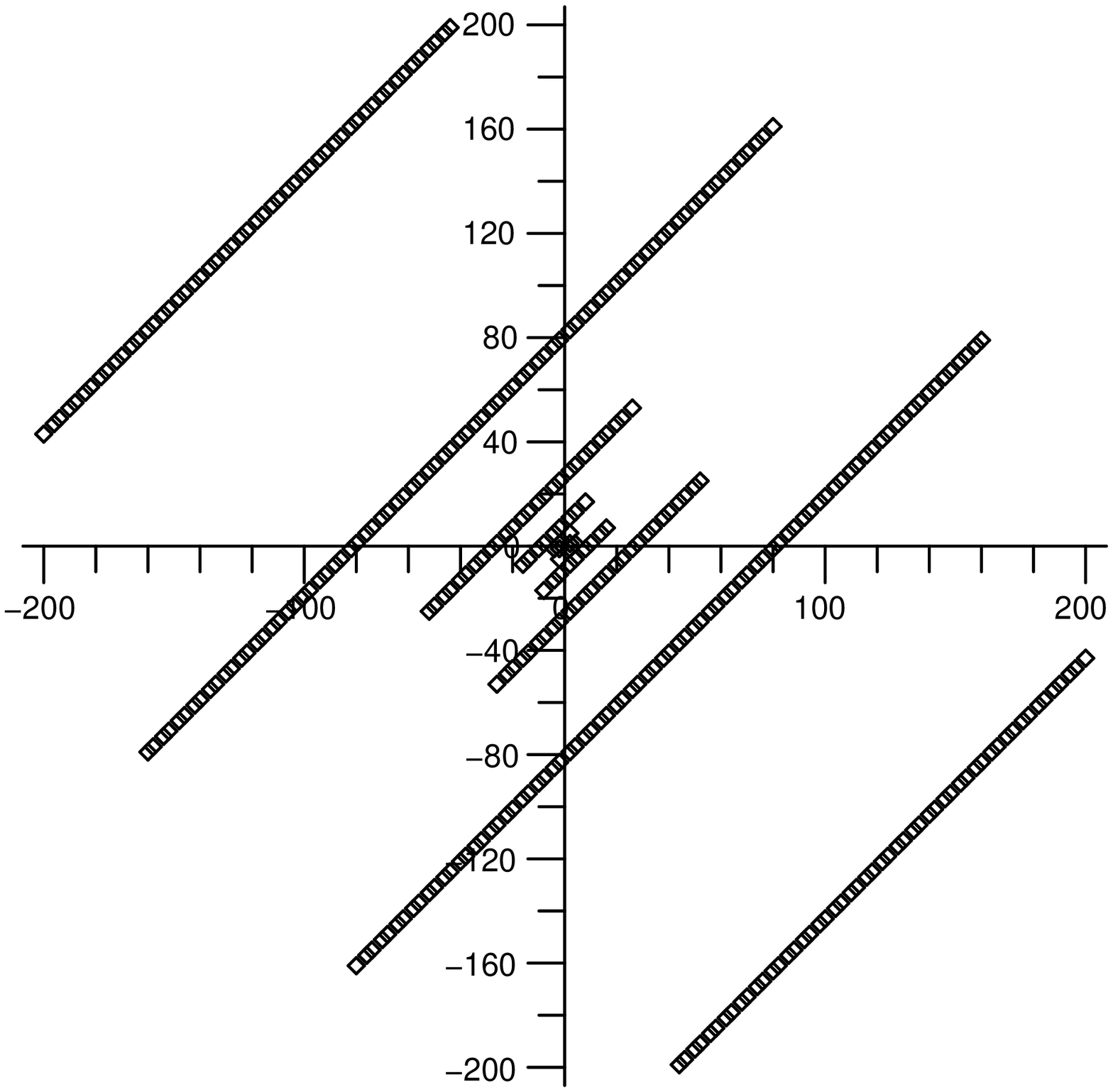,width=4in,angle=0}
  \label{FigDig}
\end{figure*}

The figure presents all valid digit sets $\{d,D\}$ for base $b=-2$ with
$-200\le d < D \le 200$. As stipulated by condition (iii) of the Theorem,
all pairs lie in one of the two obtusely angled regions bounded by $y=2x$ and
$y=\half x$.

For the proof of the Theorem we present a series of Lemmas. The first result
shows that the attractors for base $-2$ have an especially simple structure:
they are always \emph{intervals} in $\Z$.

\begin{lem} \label{Lemma3}
  Let $(\Z,-2,\{d,D\})$ be a pre-number system, with attractor ${\cal A}$, and
  suppose $d < D$. Then
  $$
    {\cal A} = \left\{ \left\lceil \frac{2d-D}3 \right\rceil, \ldots,
                  \left\lfloor \frac{2D-d}3 \right\rfloor \right\}.
  $$
\end{lem}

\begin{pf*}{Proof.}
  Theorem \ref{ThmZBd} tells us that $ \frac{2d-D}3 \le a \le \frac{2D-d}3 $
  for any $a\in{\cal A}$. We will show that these bounds are sharp. We use the
  following argument: on an arithmetic progression of difference $2$, the
  dynamic mapping $T$ is an affine linear map with slope $-\frac12$, so such a
  progression will be mapped, with its order reversed, onto an \emph{interval}.
  Thus the image of any interval $S$ under $T$ can be computed by splitting $S$
  into its even and odd parts (which are $S\cap 2\Z$ and $S\cap (2\Z+1)$,
  respectively), and considering the effect of $T$ on these parts separately.
  
  Suppose first that $d+D\equiv 0\pmod{3}$; then $2d-D$ and $2D-d$ are
  divisible by $3$. Let $a=\frac{2d-D}3$ and $A=\frac{2D-d}3$; we will prove
  that ${\cal A}=\{a,\ldots,A\}$. Note that $a\equiv D\pmod{2}$ and 
  $A\equiv d\pmod{2}$. We compute
  \begin{align*}
    T(a) &= \frac{\frac{2d-D}3-D}{-2} = \frac{2D-d}3 = A; \\
    T(a+1) &= \frac{\frac{2d-D+3}3-d}{-2} = \frac{D+d-3}6; \\
    T(A) &= \frac{\frac{2D-d}3-d}{-2} = \frac{2d-D}3 = a; \\
    T(A-1) &= \frac{\frac{2D-d-3}3-D}{-2} = \frac{D+d+3}6 = T(a+1)+1.
  \end{align*}
  It follows that the arithmetic progression $a,\,a+2,\,\ldots,\,A-1$ is mapped
  to the interval $A,\,A-1,\,\ldots,\,T(a+1)+1$, while the other progression
  $a+1,\, a+3,\,\ldots,\,A$ is mapped to $T(a+1),\,T(a+1)-1,\,\ldots,\,a$. Thus
  the interval $\{a,\ldots,A\}$ is equal to its image under $T$, which shows
  that it is equal to the attractor ${\cal A}$.
  
  Suppose that $d+D\equiv 1\pmod{3}$. Let us write $a=\lceil \frac{2d-D}3
  \rceil = \frac{2d-D+1}3$ and $A=\lfloor \frac{2D-d}3 \rfloor=\frac{2D-d-2}3$;
  we will prove that ${\cal A}=\{a, \ldots,A\}$. Note that $a\equiv D+1\equiv 
  d\pmod{2}$, and that $A\equiv d-2\equiv d\pmod{2}$. Using this, we compute
  \begin{align*}
    T(a)   &= \frac{ \frac{2d-D+1}3 - d}{-2} = \frac{D+d-1}6; \\
    T(a+1) &= \frac{ \frac{2d-D+4}3 - D}{-2} = \frac{2D-d-2}3 = A; \\
    T(A)   &= \frac{ \frac{2D-d-2}3 - d}{-2} = \frac{2d-D+1}3 = a; \\
    T(A-1) &= \frac{ \frac{2D-d-5}3 - D}{-2} = \frac{D+d+5}6 = T(a)+1.
  \end{align*}
  We again use the fact that $T$ is affine linear, with slope $-\frac12$, on
  arithmetic progressions of difference $2$. Thus, the progression
  $a,\,a+2,\,\ldots,\,A-2,\,A$ is mapped by $T$ to the interval
  $a,\,\ldots,\,T(a)$ (in reversed order), while the progression
  $a+1,\,a+3,\,\ldots,\,A-1$ is mapped to $T(a)+1,\,T(a)+2,\,\ldots,\,A$. We
  see that $\{a,\ldots,A\}$ is mapped unto itself by $T$,
  which proves the claim.

  Finally, the case where $d+D\equiv 2\pmod{3}$ is reduced to the previous by
  considering the digits $\{-d,-D\}$.
\qed \end{pf*}

\begin{lem} \label{LemCycle2}
  Let $(\Z,-2,\{d_0,d_1\})$ be a pre-number system, with attractor ${\cal A}$.
  Write $\delta = d_0 - d_1$. Then $a\in{\cal A}$ is contained in a cycle of
  length $\ell$ within ${\cal A}$ if and only if
  \begin{equation} \label{EqCycle2}
    (d_0-3a)\frac{ (-2)^\ell-1 }{-3\delta} = \sum_{i=0}^{\ell-1}
    \varepsilon_i(-2)^i
  \end{equation}
  for some $\varepsilon_i\in\{0,1\}$, and $\ell$ is minimal with this property.
\end{lem}

\begin{pf*}{Proof.}
  For any base $b$, a cycle of length $\ell$ in the attractor has the form
  \begin{displaymath}
    a_0 \rightarrow a_1=\frac{a_0-d_0}b \rightarrow a_2=
    \frac{ \frac{a_0-d_0}b - d_1 }b = \frac{a_0}{b^2} - \left( \frac{d_0}{b^2}
      + \frac{d_1}b \right) \rightarrow \ldots \rightarrow a_\ell=a_0,
  \end{displaymath}
  with $a_i\in{\cal A}$ and $d_i\in{\cal D}$ for all $i$. Continuing the
  expansion of the elements and multiplying through by $b^\ell$, we find
  \begin{displaymath}
    a_0 ( 1 - b^\ell ) = \sum_{i=0}^{\ell-1} d_ib^i.
  \end{displaymath}
  Conversely, it is clear that if $a(1-b^\ell)$ can be written in this
  form, for some $a\in{\cal A}$, and $\ell$ is minimal with this property,
  then $a$ starts a cycle of length $\ell$.

  In our case, the digits $d_i$ are either $d_0$ or $d_0 - \delta$. This gives
  \begin{displaymath}
    a_0 ( 1 - b^\ell ) = d_0 \frac{b^\ell-1}{b-1} - \delta\sum_{i=0}^{\ell-1}
    \varepsilon_i b^i,
  \end{displaymath}
  with $\varepsilon_i\in\{0,1\}$ for all $i$. It follows that
  \begin{displaymath}
    ( d_0 + (b-1)a_0 ) (b^\ell-1) = (b-1)\delta \sum_{i=0}^{\ell-1}
    \varepsilon_i b^i.
  \end{displaymath}
  The Lemma now follows by substituting $b=-2$.
\qed \end{pf*}

We will use the \emph{$q$-adic valuation} $v_q$ for a prime $q$: for an integer
$b\ne 0$, $v_q(b)$ denotes the exact number of factors $q$ in $b$. 

\begin{lem} \label{LemArtin}
  Let $q$ be an odd prime, let $b$ be an integer with $|b|\ge 2$, coprime to
  $q$, and let $n$ be a nonnegative integer. Then $q$ divides $b^n-1$ if
  and only if $\ord_q(b)$ divides $n$. If $q$ divides $b^n-1$, then
  \begin{displaymath}
    v_q(b^n-1) = v_q(n) + v_q\left( b^{\ord_q(b)}-1 \right).
  \end{displaymath}
\end{lem}

\begin{pf*}{Proof.}
This result is a special case of Lucas' \emph{law of repetition}. For a proof,
see \cite{Artin1}.
\qed \end{pf*}

\paragraph*{Example.}
Consider the digits $\{30,111\}$, so $\delta=-81$. The attractor for base $-2$
with these digits is $\{ -17,\ldots,64 \}$. Both digits are divisible by $3$,
which shows the existence of two $1$-cycles. The complete cycle structure is
\begin{align*}
  {\cal A} = & \{ 10 \} \cup 
    \{ 37 \} \cup 
    \{ -17, 64 \} \cup 
    \{ -8, 19, 46 \} \cup 
    \{ 1, 55, 28 \} \cup \mbox{{}} \\
    & \{ 4, 13, 49, 31, 40, -5, 58, -14, 22 \} \cup
    \{ -2, 16, 7, 52, -11, 61, 25, 43, 34 \} \mbox{{}} \\ 
    & \{ 0, 15, 48, -9, 60, -15, 63, 24, 3, 54, -12, 21, 45, 33, 39, 36, -3,
    57, 27, 42, \\
    & \quad -6, 18, 6, 12, 9, 51, 30 \} \cup \mbox{{}} \\
    & \{ -1, 56, -13, 62, -16, 23, 44, -7, 59, 26, 2, 14, 8, 11, 50, -10, 20,
    5, 53, 29, \\
    & \quad 41, 35, 38, -4, 17, 47, 32 \}.
\end{align*}
Of these, the cycle lengths $\ell$ that are powers of $3$ are not that
surprising, because $(-2)^\ell-1$ is then divisible by $\ell$, and the
remaining factors of the denominator $3\delta$ are found in $(d_0-3a)$. The
$2$-cycle is legitimised by the following calculation: the factor $(-2)^2-1$
cancels the $3$ in the denominator, while we have $d_0-3\cdot(-17)=81$ and
$d_0-3\cdot 64=-162$, both of which are divisible by $\delta$.

\begin{lem} \label{LemOneCycle}
  Let $(\Z,-2,\{d_0,d_1\})$ be a pre-number system, with attractor ${\cal A}$.
  Then ${\cal A}$ consists of exactly one cycle if and only if either
  $|d_0-d_1|=1$, or
  \begin{romanlist}
    \item $|d_0-d_1|=3^i$ for some $i\ge 1$, and
    \item $3\nmid d_0$ and $3\nmid d_1$.
  \end{romanlist}
\end{lem}

\begin{pf*}{Proof.}
  Write $d_1 = d_0 - \delta$ as above, so $\delta$ is an odd integer. We first
  prove the ``if''-part.

  First, assume $|\delta|=1$. If $d_0+d_1\not\equiv 0\pmod{3}$, then by Lemma
  \ref{Lemma3}, ${\cal A}$ consists of only one element, and the claim is
  obvious. If $d_0+d_1\equiv 0\pmod{3}$, then ${\cal A}$ has $2$ elements,
  again by Lemma \ref{Lemma3}. If the claim fails, there must be a $1$-cycle in
  ${\cal A}$, and this implies that either $d_0$ or $d_1$ is divisible by $3$,
  by Corollary \ref{CorOneCycle}. But this contradicts the assumption that
  $d_0+d_1\equiv 0\pmod{3}$. It follows that ${\cal A}$ has a single $2$-cycle,
  as desired.

  Next, assume $|\delta|=3^i$ for some $i\ge 1$. As remarked earlier, if $3$
  divides either $d_0$ or $d_1$, we immediately obtain a $1$-cycle in 
  ${\cal A}$. Therefore we exclude this case, and it follows that
  $d_0+d_1\not\equiv 0\pmod{3}$. By Lemma \ref{Lemma3}, we conclude that
  \begin{displaymath}
    |{\cal A}| = |\delta| = 3^i.
  \end{displaymath}
  Let $\ell$ be the length of the longest cycle in ${\cal A}$. By Lemma
  \ref{LemCycle2}, and because $3\nmid d_0$, we conclude that $3^{i+1} \mid
  (-2)^\ell - 1$. Now by Lemma \ref{LemArtin}, taking $b=-2$ and $q=3$, we find
  that
  \begin{displaymath}
    3^{i+1} \mid (-2)^\ell-1 \Rightarrow 3^i \mid \ell.
  \end{displaymath}
  Because $\ell\le|\delta|$, it follows that $\ell=|\delta|$, so that ${\cal
  A}$ consists of just one cycle, and the first half of the Lemma is proved.

  ~

  Now we prove the ``only if''-part. Suppose that ${\cal A}$ consists of just
  one cycle. We distinguish two cases, namely whether $3$ divides $\delta$ or
  not.

  First, assume that $3$ divides $\delta$. Now either both $d_0$ and $d_1$ are
  divisible by $3$, or neither of them is. If both are divisible by $3$, then
  the attractor has two distinct $1$-cycles, which is a contradiction. Thus,
  $3$ divides neither of $d_0$ and $d_1$. By Lemma \ref{Lemma3}, we find that
  ${\cal A}$ is an interval of length $|\delta|$, so that we have just one
  cycle of length $|\delta|$.

  Now consider \eqref{EqCycle2}. Because ${\cal A}$ contains an element from 
  every residue class modulo $\delta$, and because $3\nmid d_0$, we can choose
  $a_0\in{\cal A}$ so that $\gcd(d_0-3a_0,\delta)=1$. It follows that 
  \begin{displaymath}
    3\delta \mid (-2)^{|\delta|}-1,
  \end{displaymath}
  and this does not hold for any smaller exponent than $|\delta|$. We will show
  that this implies that $|\delta|$ is a power of $3$.

  The assumption means that the order of $-2$ in the multiplicative group
  $\left(\Z/3\delta\Z\right)^*$ is equal to $|\delta|$. But this order divides
  the order of the group, which is $\phi(3|\delta|)=3\phi(|\delta|)$, as we
  assume that $3\mid\delta$. Let $p$ be the largest prime divisor of $\delta$,
  and suppose $p>3$. Then $\phi(3|\delta|)$ has less factors $p$ than $\delta$,
  so that the divisibility relation is impossible. It follows that $\delta$ is
  a power of $3$.
  
  Finally, assume that $3$ does not divide $\delta$. If $3$ divides $d_0$, then
  $3$ does not divide $d_1$, and ${\cal A}$ has exactly one $1$-cycle. It
  follows that ${\cal A}$ has just one element. Also, we have
  $d_0+d_1\not\equiv 0\pmod{3}$, so $|{\cal A}|=|\delta|$ by Lemma
  \ref{Lemma3}. We obtain $|\delta|=1$, as desired. 
  
  If $3$ divides neither of $d_0$ or $d_1$, then one easily verifies that
  $d_0+d_1\equiv 0\pmod{3}$. In this case, Lemma \ref{Lemma3} shows that
  $\frac{2D-d}3$ and $\frac{2d-D}3$ are in ${\cal A}$. But these two elements
  constitute a $2$-cycle under $T$, and it follows that ${\cal A}$ has just
  these two elements. As $|{\cal A}|$ is equal to $|\delta|+1$, again by Lemma
  \ref{Lemma3}, we see that $|\delta|=1$, as desired.
\qed \end{pf*}

\begin{pf*}{Proof of Theorem \ref{ThmMinus2}}
The condition of having one even and one odd digit is obviously necessary. Now
the number system condition is equivalent to the requirement that the attractor
${\cal A}$ consists of exactly one cycle under the dynamic map $T$, and that
this cycle contains $0$.

By Lemma \ref{LemOneCycle}, the attractor has one cycle if and only if $D-d=1$,
or $D-d=3^i$ for some $i\ge 1$ and neither $D$ nor $d$ is divisible by $3$.
Next, Lemma \ref{Lemma3} tells us whether $0$ is in the attractor, as follows.

If $D-d=1$ and $3$ divides one of the digits, we have $D+d\not\equiv 0\pmod{3}$,
so ${\cal A}$ consists of just one element. If $3\mid D$, then this element is
$-D/3$, and if $3\mid d$, it is $-d/3$, as these elements generate $1$-cycles.
It follows that the digit divisible by $3$ must be $0$.

If $D-d=1$ and $3$ does not divide a digit, then $D+d\equiv 0\pmod{3}$, so 
${\cal A}$ has just the elements $\frac{2d-D}3$ and $\frac{2D-d}3$, forming a
$2$-cycle. One of these elements is $0$, and one verifies that $2d\le D$ and
$2D\ge d$ are necessary and sufficient conditions for this to hold.

If $D-d=3^i$ for $i\ge 1$, and $3$ does not divide a digit, then $D+d\not\equiv
0\pmod{3}$. Here again, from the form of ${\cal A}$ given by Lemma
\ref{Lemma3}, one easily verifies that the two conditions $2d\le D$ and $2D\ge
d$ exactly ensure that $0\in{\cal A}$.
\qed \end{pf*}

\begin{ack}
The research leading to this paper was supported by the Austrian Science
Foundation FWF, project S9606, which is part of the Austrian National Research
Network ``Analytic Combinatorics and Probabilistic Number Theory.''

The anonymous referee is to be thanked for his/her useful suggestions.
\end{ack}

\hyphenation{ma-the-ma-ti-ques}\def\cprime{$'$}
  \def\polhk#1{\setbox0=\hbox{#1}{\ooalign{\hidewidth
  \lower1.5ex\hbox{`}\hidewidth\crcr\unhbox0}}} \def\cprime{$'$}

\end{document}